# The central limit theorem under random truncation

WINFRIED STUTE[1] and JANE-LING WANG[2]

[1]*Mathematical Institute, University of Giessen, Arndtstr. 2, D-35392 Giessen, Germany.
E-mail: winfried.stute@math.uni-giessen.de*
[2]*Department of Statistics, University of California, Davis, CA 95616, USA.*

Under left truncation, data $(X_i, Y_i)$ are observed only when $Y_i \leq X_i$. Usually, the distribution function $F$ of the $X_i$ is the target of interest. In this paper, we study linear functionals $\int \varphi \, dF_n$ of the nonparametric maximum likelihood estimator (MLE) of $F$, the Lynden-Bell estimator $F_n$. A useful representation of $\int \varphi \, dF_n$ is derived which yields asymptotic normality under optimal moment conditions on the score function $\varphi$. No continuity assumption on $F$ is required. As a by-product, we obtain the distributional convergence of the Lynden-Bell empirical process on the whole real line.

## 1. Introduction and main results

In this paper, we provide some further methodology for statistical analysis of data which are truncated from the left. To be more specific, let $(X_i, Y_i), 1 \leq i \leq N$, denote a sample of independent bivariate data such that, for each $i$, $X_i$ is also independent of $Y_i$. Denote by $F$ and $G$, respectively, the unknown distribution functions of $X$ and $Y$. Typically, $F$ is the target of interest. Now, under left truncation, $X_i$ is observed only when $Y_i \leq X_i$. As a consequence, the empirical distribution of the $X$'s is unavailable and cannot serve as a basic process to compute other statistics.

The nonparametric maximum likelihood estimator of $F$ for left-truncated data was first derived by Lynden-Bell (1971). Its first mathematical investigation may be attributed to Woodroofe (1985), who also reviewed some examples of truncated data from astronomy and economics; see also Wang (1989) for applications in the analysis of AIDS data.

Now, denoting by $n$ the number of data which are actually observed, we have, by the strong law of large numbers (SLLN),

$$\frac{n}{N} \to \alpha \equiv \mathbb{P}(Y \leq X) \qquad \text{as } N \to \infty \text{ with probability one.}$$

Without further mention, we shall assume that $\alpha > 0$ because, otherwise, nothing will be observed. Of course, $\alpha$ will be unknown. Conditionally on $n$, the observed data are still independent, but the joint distribution of $X_i$ and $Y_i$ becomes

$$H^*(x,y) = \mathbb{P}(X \leq x, Y \leq y | Y \leq X) = \alpha^{-1} \int_{(-\infty,x]} G(y \wedge z) F(\mathrm{d}z).$$

The marginal distribution of the observed $X$'s thus equals

$$F^*(x) \equiv \alpha^{-1} \int_{(-\infty,x]} G(z) F(\mathrm{d}z). \tag{1.1}$$

It may be consistently estimated by the empirical distribution function of the known $X$'s:

$$F_n^*(x) = \frac{1}{n} \sum_{i=1}^n \mathbf{1}_{\{X_i \leq x\}}, \qquad x \in \mathbb{R}.$$

The problem, however, is one of reconstructing $F$ and not $F^*$ from the available data $(X_i, Y_i), 1 \leq i \leq n$. A crucial quantity in this context is the function

$$C(z) = \mathbb{P}(Y \leq z \leq X | Y \leq X) = \alpha^{-1} G(z) [1 - F(z-)], \tag{1.2}$$

where

$$F(z-) = \lim_{x \uparrow z} F(x)$$

is the left-continuous version of $F$ and $F\{z\} = F(z) - F(z-)$ is the $F$-mass at $z$. The function $C$ may be consistently estimated by

$$C_n(z) = n^{-1} \sum_{i=1}^n \mathbf{1}_{\{Y_i \leq z \leq X_i\}}.$$

It is very helpful to express the cumulative hazard function of $F$,

$$\Lambda(x) = \int_{(-\infty,x]} \frac{F(\mathrm{d}z)}{1 - F(z-)},$$

in terms of estimable quantities. For this, let

$$a_G = \inf\{x : G(x) > 0\}$$

be the largest lower bound for the support of $G$. Similarly for $F$. From (1.1) and (1.2), we obtain

$$\int_{(a_G,x]} \frac{F(\mathrm{d}z)}{1 - F(z-)} = \int_{(a_G,x]} \frac{F^*(\mathrm{d}z)}{C(z)}. \tag{1.3}$$

Provided that $a_G \leq a_F$ and $F\{a_F\} = 0$, the left-hand side equals $\Lambda(x)$. Otherwise, this is no longer true and, as Woodroofe (1985) pointed out, $F$ cannot be fully recovered from the available data. The situation is similar for right-censored data; see Stute and Wang (1993) for a detailed discussion for (upper) boundary effects there.

Throughout the paper, we shall therefore assume that

$$a_G \leq a_F \quad \text{and} \quad F\{a_F\} = 0. \tag{1.4}$$

If $a_G < a_F$, the second assumption is superfluous, while it is automatically satisfied when $F$ is continuous. For the moment, however, no other assumptions such as continuity of $F$ or $G$ will be needed.

In differential terms, equation (1.3) leads to

$$F(\mathrm{d}x) = (1 - F(x-))\frac{F^*(\mathrm{d}x)}{C(x)}. \tag{1.5}$$

The Lynden-Bell (1971) estimator $F_n$ of $F$ is obtained as the solution of the so-called self-consistency equation, that is, as the solution of (1.5) after having replaced $F^*$ and $C$ with $F_n^*$ and $C_n$, respectively:

$$F_n(\mathrm{d}x) = (1 - F_n(x-))\frac{F_n^*(\mathrm{d}x)}{C_n(x)}. \tag{1.6}$$

Solving for $F_n$ yields the product integral representation of $F_n$:

$$1 - F_n(x) = \prod_{\text{distinct } X_i \leq x} \left[1 - \frac{F_n^*\{X_i\}}{C_n(X_i)}\right]. \tag{1.7}$$

If there are no ties among the $X$'s, (1.7) simplifies to become

$$1 - F_n(x) = \prod_{X_i \leq x} \left[\frac{nC_n(X_i) - 1}{nC_n(X_i)}\right]. \tag{1.8}$$

Note that $nC_n(X_i) \geq 1$, so each ratio is well defined. Since, in this paper, our objective will be to study general linear statistics based on $F_n$, namely Lynden-Bell integrals $\int \varphi \, \mathrm{d}F_n$, we also introduce the Lynden-Bell weights $W_{in}$ attached to each datum in the $X$-sample. For this, denote by $X_{1:n} < \cdots < X_{m:n}$ the $m$ distinct order statistics, with $m$ possibly strictly less than $n$.

From (1.7), we obtain

$$W_{in} \equiv F_n\{X_{i:n}\} = [1 - F_n(X_{i-1:n})]\frac{F_n^*\{X_{i:n}\}}{C_n(X_{i:n})} \tag{1.9}$$

and

$$\int \varphi \, \mathrm{d}F_n = \sum_{i=1}^m W_{in} \varphi(X_{i:n}). \tag{1.10}$$

For $\varphi = 1_{(-\infty,x]}$, we are back at $F_n(x)$. As for other $\varphi$'s, we refer to Stute and Wang (1993) or Stute (2004), who considered possible applications of empirical integrals in the context of right-censored data. More general statistical functionals often admit expansions, in which the leading term is of the form $\int \varphi \, dF_n$, with $\varphi$ denoting the associated influence function. Since, for a fully observable data set with $F_n^e$ being the classical empirical distribution function, $\int \varphi \, dF_n^e$ is just a sample mean to which, under a second moment assumption, the central limit theorem (CLT) applies, distributional convergence of $\int \varphi \, dF_n$ therefore constitutes an extension of the CLT to the left-truncation case. The corresponding SLLN has been studied in papers by He and Yang (1998a, 1998b). The CLT for censored data is due to Stute (1995).

In the present situation, the CLT is much more elusive than for randomly censored data. In the censored data case, the right tails create technical difficulties, but not the left. For left truncation, however, both sides create problems. This is already seen via the functions $C$ and $C_n$ which decrease to zero on the left and on the right tail, and with $C_n$ appearing in denominators. The only trivial bound is $nC_n(X_i) \geq 1$, which keeps everything from being "not well defined". Moreover, $C$ and $C_n$ are not monotone, again contrary to the censored data case, where the role of the $C$'s is played by survival functions. This non-monotone feature of $C_n$ creates additional technical complications for truncated data. Worse than that, $C_n$ may also become zero between the data points in the central part. Keiding and Gill (1990), who recognized the danger of these sets, have called these holes the "empty inner risk sets". Not all authors seem to know about these problems since a detailed study of the $C_n$ process is sometimes missing. A consequence of these empty inner risk sets, as revealed by the representation (1.7), is the loss of mass of $F_n$ after the first such hole. More specifically, suppose that, in terms of Keiding and Gill (1990), a risky hole exists at $X_j$. This means that $X_j$ is not covered by any other pair $(X_i, Y_i)$. Hence, $nC_n(X_j) = nF_n^*\{X_j\}$. From (1.7), we get that all data points right to $X_j$ have mass zero under $F_n$. In proofs, this disallows the incorporation of exponential representations of the weights. To circumvent this difficulty, we first study an asymptotically equivalent estimator $\hat{F}_n$. This estimator is defined via

$$\int \varphi \, d\hat{F}_n = \sum_{i=1}^{m} \frac{\varphi(X_{i:n}) F_n^*\{X_{i:n}\}}{C_n(X_{i:n})} \prod_{j=1}^{i-1} \left[ 1 - \frac{nF_n^*\{X_{j:n}\}}{nC_n(X_{j:n}) + 1} \right]. \tag{1.11}$$

The extra summand 1 in the denominator allows for contributions from data which have holes on their left. As we shall see in a small simulation study, it may have a robustifying effect, resulting in smaller mean squared error (MSE) for small sample sizes. Asymptotically, $\int \varphi \, dF_n$ and $\int \varphi \, d\hat{F}_n$ are equivalent at the $n^{1/2}$-rate, which facilitates the asymptotic theory for $\int \varphi \, dF_n$. Another technical problem is caused by possible ties. In such a situation, (1.8) does not apply. Lemma 1.1 will show how the general case covered by (1.7) and (1.9) may be traced back to the case of a continuous $F^*$.

The main result of this paper, Theorem 1.1, provides a representation of $\int \varphi \, dF_n$ as a sum of i.i.d. random variables under minimal assumptions on $\varphi$ and the truncation mechanism. Usually, linear i.i.d. representations of complicated estimators will include the Hájek projection of the statistic of interest. In the case of (1.10), however, this projection can be computed only up to remainders, and this is exactly what Theorem 1.1 does. More precisely, the proof of Theorem 1.1 proceeds in two steps. In the first step we have to identify all error terms which are negligible. The leading terms will turn out to be V-statistics; see Serfling (1980). Finally, an application of the Hájek projection to the leading terms yields the desired i.i.d. representation. Needless to say, asymptotic normality follows immediately. We shall also add some interesting comments on a so-called uniform representation. Proofs will be given in Section 3.

Theorem 1.1 will hold under the following two assumptions:

(A)    (i)  $\int \varphi^2/G \, dF < \infty$;
       (ii)  $\int \frac{dF}{G} < \infty$.

Condition (ii) already appeared in Woodroofe (1985) in his study of the Lynden-Bell process, that is, when he considered indicators $\varphi = 1_{(-\infty, x_0]}$. Non-technically speaking, it is needed to guarantee enough information in the left tails to estimate $F$ at the rate $n^{-1/2}$. Under slightly stronger assumptions, Stute (1993) obtained an almost sure representation with sharp bounds on the remainder, again for indicators; see also Chao and Lo (1988). Condition (i) guarantees, among other things, that the leading terms in the i.i.d. representation admit a finite second moment so that asymptotic normality holds. Since $G \leq 1$, it implies $\int \varphi^2 \, dF < \infty$, which is the standard finite moment assumption when no truncation occurs. When $\varphi$ has a finite second $F$-moment and is locally bounded in a neighborhood of $a_G$, then (i) is implied by (ii). Note also that (i) and (ii) are always satisfied when $a_G < a_F$ and $\int \varphi^2 \, dF < \infty$. A CLT for truncated data is also contained in Sánchez Sellero et al. (2005). Apart from continuity assumptions, they also need conditions which, in our notation, require finiteness of the integral

$$\int \varphi_0^2(x)(1 - F(x))^{-5} F(dx),$$

where $|\varphi| \leq \varphi_0$. Since, however, this integral equals infinity for constant $\varphi_0$, their result is not applicable to bounded $\varphi$'s, not to mention $\varphi$'s which increase to infinity as $x \to \infty$. Rather, finiteness of the above integral is only obtained for $\varphi_0$'s which converge to zero fast enough in the right tails.

The focus of this paper is, however, on distributional convergence for which (A) will suffice. Theorem 1.1 is formulated for the case when $F$ is continuous. This guarantees that among the observed $X$'s, there will be no ties, with probability one. In such a situation, we obtain

$$\int \varphi \, d\hat{F}_n = \int \frac{\varphi(x)\gamma_n(x)}{C_n(x)} F_n^*(dx), \tag{1.12}$$

with

$$\gamma_n(x) = \exp\left\{ n \int_{-\infty}^{x-} \ln\left[1 - \frac{1}{nC_n(y) + 1}\right] F_n^*(dy) \right\}. \tag{1.13}$$

At the end of this section, we shall show how general Lynden-Bell integrals may be traced back to the present case.

**Theorem 1.1.** *Under Assumptions* (A) *and (1.4), assume that $F$ is continuous. We then have*

$$\int \varphi[dF_n - dF] = \int \frac{\psi(y)}{C(y)}[F_n^*(dy) - F^*(dy)]$$
$$- \int \frac{C_n(y) - C(y)}{C^2(y)}\psi(y)F^*(dy) + o_{\mathbb{P}}(n^{-1/2}),$$

*where*

$$\psi(y) = \varphi(y)(1 - F(y)) - \int_{\{y<x\}} \frac{\varphi(x)(1-F(x))}{C(x)} F^*(dx) = \int_{\{y<x\}} [\varphi(y) - \varphi(x)]F(dx).$$

For indicators $\varphi = 1_{(-\infty, x_0]}$, the leading term already appears in Theorem 2 in Stute (1993).

***Remark 1.1.*** If one checks the proof of Theorem 1.1 step by step, the following fact will be revealed. If, rather than a single $\varphi$, one considers a collection $\{\varphi\}$, then the error terms are uniformly small in the sense that the remainder is $o_{\mathbb{P}}(n^{-1/2})$ uniformly in $\varphi$, provided that $|\varphi| \leq \varphi_0$ for some $\varphi_0$ satisfying $\int \frac{\varphi_0^2}{G} dF < \infty$. Actually, all remainders may be bounded from above in absolute value by replacing $|\varphi|$ by $\varphi_0$. Compared with Sánchez Sellero *et al.* (2005), no VC property for the $\varphi$'s is needed for a representation as a V-statistic process. For the i.i.d. representation, one has to guarantee that the errors in the Hájek projection are also uniformly small. These errors, however, form a class of degenerate U-statistics. This kind of process was studied in Stute (1994) and the achieved bounds are useful for handling the error terms in the second half of the proof. Details are omitted.

**Corollary 1.1.** *Under the assumptions of Theorem 1.1, we have*

$$n^{1/2} \int \varphi[dF_n - dF] \to \mathcal{N}(0, \sigma^2) \text{ in distribution,}$$

*with*

$$\sigma^2 = \text{Var}\left\{\frac{\psi(X)}{C(X)} - \int_Y^X \frac{\psi(y)}{C^2(y)} F^*(dy)\right\}.$$

It is not difficult to see that $\sigma^2 < \infty$ under (A).

Finally, if, in Remark 1.1, we take for $\{\varphi\}$ the class of all indicators $\varphi = 1_{(-\infty, x]}$ and set $\varphi_0 \equiv 1$, we obtain the following corollary.

**Corollary 1.2.** *Under $\int \frac{dF}{G} < \infty$ and (1.4), and for a continuous $F$, we have*

$$F_n(x) - F(x) = \int \frac{\psi_x}{C}[dF_n^* - dF^*]$$

$$-\int \frac{C_n - C}{C^2} \psi_x F^*(\mathrm{d}y) + o_{\mathbb{P}}(n^{-1/2})$$

*uniformly in $x$, where $\psi_x$ is the $\psi$ belonging to $\varphi = 1_{(-\infty, x]}$.*

**Remark 1.2.** To make the point clear, Corollary 1.2 provides a representation which holds uniformly on the whole real line and not only on subintervals $(-\infty, b]$ with $b < b_F = \sup\{x : F(x) < 1\}$, as is usually the case in the literature.

A major technical problem for proving Theorem 1.1 for a general $F$ is caused by the fact that for discontinuous $F$, ties may arise. As before, denote by $X_{1:n} < \cdots < X_{m:n}$ the $m$ ordered distinct data in the observed $X$-sample. To circumvent ties, one may use the fact that each $X_i$ may be written in the form $X_i = F^{*-1}(U_i)$, where $U_1, \ldots, U_n$ is a sample of independent random variables with a uniform distribution on $(0, 1)$ and $F^{*-1}$ is the quantile function of $F^*$. The construction of the $U$'s is similar to the construction in Lemma 2.8 of Stute and Wang (1993), with $H$ there replaced by $F^*$ here.

For the following, recall that a quantile function is continuous from the left. Furthermore, with probability one, $F^{*-1}$ is also right-continuous at each $U_i$. With this in mind, one can see that the corresponding truncating sample for the $U_i$'s consists of $F^*(Y_i-), 1 \leq i \leq n$. If we denote by $C_n^U$ the $C_n$-function corresponding to the pseudo-observations $(U_i, F^*(Y_i-)), 1 \leq i \leq n$, and if we let $U_{i1:n} < \cdots < U_{id_i:n}$ denote the ordered $U$'s that satisfy $F^{*-1}(U_{ij:n}) = X_{i:n}$, with $d_i = nF_n^*\{X_{i:n}\}$, then

$$nC_n^U(U_{i1:n}) = nC_n(X_{i:n}) \tag{1.14}$$

and

$$nC_n^U(U_{ij:n}) = nC_n^U(U_{i,j-1:n}) - 1 \qquad \text{for } 2 \leq j \leq d_i. \tag{1.15}$$

Also, note that the Lynden-Bell estimator $F_n^U$ for the pseudo-observations satisfies (1.8), that is,

$$1 - F_n^U(u) = \prod_{U_i \leq u} \left[\frac{nC_n^U(U_i) - 1}{nC_n^U(U_i)}\right].$$

Introducing the function $\varphi^* = \varphi \circ F^{*-1}$, the analog of (1.10) becomes

$$\int \varphi^* \, \mathrm{d}F_n^U = \sum_{i=1}^m \sum_{j=1}^{d_i} \varphi(F^{*-1}(U_{ij:n})) \frac{1 - F_n^U(U_{ij:n}-)}{nC_n^U(U_{ij:n})}$$

$$= \sum_{i=1}^m \varphi(X_{i:n}) \sum_{j=1}^{d_i} \frac{1 - F_n^U(U_{ij:n}-)}{nC_n^U(U_{ij:n})}.$$

In Lemma 1.1, we will show that for each $1 \leq i \leq m$ and every $1 \leq j \leq d_i$, we have

$$\frac{1 - F_n^U(U_{ij:n}-)}{nC_n^U(U_{ij:n})} = \frac{1 - F_n(X_{i-1:n})}{nC_n(X_{i:n})}, \tag{1.16}$$

where $X_{0:n} = -\infty$. It follows from (1.9), (1.10) and (1.16) that

$$\int \varphi^* \, \mathrm{d}F_n^U = \int \varphi \, \mathrm{d}F_n \qquad \text{with probability one.}$$

In conclusion, the study of Lynden-Bell integrals may be traced back to the case when the variables of interest are uniform on $(0,1)$ and, therefore, with probability one, have no ties. At the same time, our handling of ties does not require external randomization, but maintains the product limit structure of the weights and hence of the estimators.

**Lemma 1.1.** *For each $1 \leq i \leq m$ and $1 \leq j \leq d_i$, equation (1.16) holds.*

## 2. Simulations

It is interesting to compare the small sample size behaviors of $F_n$ and $\hat{F}_n$. In a simulation study, we considered $\varphi(x) = x$, that is, the target was the mean lifetime of $X$. This $\varphi$ is the canonical score function in the classical CLT. Recall that via truncation, there is a sampling bias which would result in an upward bias if we were to take the empirical distribution function of the $X$'s and not $F_n$ or $\hat{F}_n$. Introducing more or less complicated weights has the effect, among other things, that compared with the empirical distribution function, the bias is reduced.

In the following, we report on some simulation results which are part of a much more extensive study. Typically, for this $\varphi$, $\int \varphi \, \mathrm{d}\hat{F}_n$ outperforms $\int \varphi \, \mathrm{d}F_n$ in terms of the MSE when $10 \leq n \leq 40$ and truncation is heavy. For larger $n$ ($n \geq 40$), the difference is negligible. MSE was computed via Monte Carlo. See Table 1 for details.

In the simulations, both $X$ and $Y$ were exponentially distributed with parameter 1.

**Table 1.** Comparison of MSE

| Sample size | MSE($\int \varphi \, \mathrm{d}F_n$) | MSE($\int \varphi \, \mathrm{d}\hat{F}_n$) |
|---|---|---|
| $n = 10$ | 0.93 | 0.56 |
| $n = 20$ | 0.45 | 0.39 |
| $n = 30$ | 0.38 | 0.32 |
| $n = 40$ | 0.35 | 0.31 |
| $n = 50$ | 0.29 | 0.27 |
| $n = 60$ | 0.26 | 0.24 |
| $n = 70$ | 0.24 | 0.23 |
| $n = 80$ | 0.22 | 0.22 |
| $n = 90$ | 0.22 | 0.21 |
| $n = 100$ | 0.20 | 0.19 |

## 3. Proofs

To prove Theorem 1.1, note that under continuity of $F$, (1.8) applies. Also, we may assume, without loss of generality, that all data are non-negative. Hence, all integrals appearing hereafter will be over the positive real line.

In our first lemma, we provide a bound for the function $C/C_n$. This will be needed in proofs to handle negligible terms. Although, by Glivenko–Cantelli, $C_n - C \to 0$ uniformly, bounding the ratio is more delicate in view of possible holes and the non-monotonicity of $C$ and $C_n$.

**Lemma 3.1.** *Assume $F$ is continuous. For any $\lambda$ such that $\alpha\lambda \geq 1$, one has*

$$\mathbb{P}\left(\sup_{X_i: X_i \geq b} \frac{C(X_i)}{C_n(X_i)} \geq \lambda\right) \leq \lambda e^2 \exp[-G(b)\alpha\lambda].$$

**Proof.** The proof is similar to that of Lemma 1.2 in Stute (1993), which provides an exponential bound for the supremum extended over the left tails $X_i \leq b$ rather than the right tails. □

Together with the aforementioned bound from Stute (1993), Lemma 3.1 immediately implies

$$\sup_{1 \leq i \leq n} \frac{C(X_i)}{C_n(X_i)} = O_\mathbb{P}(1) \qquad \text{as } n \to \infty. \tag{3.1}$$

Assertion (3.1) will be of some importance in forthcoming proofs as it will allow us to replace the random $C_n$ appearing in denominators with the deterministic $C$.

We are now ready to expand $\int \varphi \, d\hat{F}_n$. By (1.12) and (1.13),

$$\int \varphi \, d\hat{F}_n = \int_0^\infty \frac{\varphi(x)}{C_n(x)} \exp\left\{n \int_0^{x-} \ln\left[1 - \frac{1}{nC_n(y) + 1}\right] F_n^*(dy)\right\} F_n^*(dx)$$

$$\equiv \int_0^\infty \frac{\varphi(x)}{C_n(x)} \gamma_n(x) F_n^*(dx).$$

Set

$$\gamma(x) = 1 - F(x) = \exp\left\{-\int_0^x \frac{F^*(dy)}{C(y)}\right\}.$$

From Taylor's expansion, we obtain

$$\gamma_n(x) = \gamma(x) + e^{\Delta_n(x)}[B_n(x) + D_{n1}(x) + D_{n2}(x)],$$

where

$$B_n(x) = n \int_0^{x-} \ln\left[1 - \frac{1}{nC_n(y) + 1}\right] F_n^*(dy) + \int_0^{x-} \frac{F_n^*(dy)}{C_n(y) + 1/n},$$

$$D_{n1}(x) = -\int_0^{x-} \frac{F_n^*(\mathrm{d}y) - F^*(\mathrm{d}y)}{C(y)},$$

$$D_{n2}(x) = \int_0^{x-} \frac{[C_n(y) + 1/n - C(y)]}{[C_n(y) + 1/n]C(y)} F_n^*(\mathrm{d}y)$$

and $\Delta_n(x)$ is between the exponents of $\gamma_n(x)$ and $\gamma(x)$. Particularly, we have $\Delta_n(x) \leq 0$. Setting

$$S_{n1} = \int_0^\infty \frac{\varphi(x)}{C_n(x)} e^{\Delta_n(x)} B_n(x) F_n^*(\mathrm{d}x),$$

$$S_{n2} = \int_0^\infty \frac{\varphi(x)\gamma(x)}{C_n(x)} F_n^*(\mathrm{d}x),$$

$$S_{n3} = \int_0^\infty \frac{\varphi(x)}{C_n(x)} e^{\Delta_n(x)} D_{n1}(x) F_n^*(\mathrm{d}x)$$

and

$$S_{n4} = \int_0^\infty \frac{\varphi(x)}{C_n(x)} e^{\Delta_n(x)} D_{n2}(x) F_n^*(\mathrm{d}x),$$

we thus get

$$\int \varphi \, \mathrm{d}\hat{F}_n = S_{n1} + S_{n2} + S_{n3} + S_{n4}. \tag{3.2}$$

In the next lemma, we study the functions $D_{n1}$ and $D_{n2}$ more closely. To motivate the following, note that for each fixed $x_0$ such that $F(x_0) < 1$,

$$D_{n1}(x_0) \to 0 \quad \text{with probability one.}$$

Actually, by standard Glivenko–Cantelli arguments,

$$D_{n1}(x) \to 0 \quad \text{with probability one uniformly on } x \leq x_0.$$

Similarly, when we consider the standardized processes $x \to n^{1/2} D_{n1}(x)$, it is easy to show their distributional convergence in the Skorokhod space $D[0, x_0]$. Things change, however, if we study $D_{n1}$ on the whole support of $F^*$. Since

$$\int_0^\infty \frac{F^*(\mathrm{d}y)}{C(y)} = \infty,$$

we cannot expect uniform convergence on the whole support of $F^*$. The situation is similar for the cumulative hazard function $\Lambda$, where uniform convergence of the Nelson–Aalen estimator $\Lambda_n$ may be obtained only on compact subsets of the support of $F$.

When one evaluates these processes at $x = X_i$, though, it is known that the uniform deviation between $\Lambda_n$ and $\Lambda$ does not got to zero, but remains at least bounded; see

Theorem 2.1 in Zhou (1991). Similar things turn out to be true for $D_{n1}$ and $D_{n2}$, as Lemma 3.2 will show. Our proofs are different, though, since compared with Zhou (1991), we shall apply a truncation technique which in proofs guarantees that the suprema of $D_{n1}$ and $D_{n2}$ are bounded on large, but not too large, sets.

**Lemma 3.2.** *We have*

$$\sup_{1 \leq i \leq n} |D_{n1}(X_i)| = O_\mathbb{P}(1) \tag{3.3}$$

*and*

$$\sup_{1 \leq i \leq n} |D_{n2}(X_i)| = O_\mathbb{P}(1). \tag{3.4}$$

**Proof.** Assume, without loss of generality, that $F^*$ has unbounded support on the right. Otherwise, replace $\infty$ by $b_F^* = \sup\{x : F^*(x) < 1\}$. For a given $\varepsilon > 0$, one may find some small $c = c_\varepsilon$ and a sequence $a_n \to \infty$ such that $1 - F^*(a_n) = \frac{c}{n}$ and $\mathbb{P}(X_{n:n} \leq a_n) \geq 1 - \varepsilon$. Actually, this follows from

$$\mathbb{P}(X_{n:n} \leq a_n) = \left[1 - \frac{c}{n}\right]^n \to \exp(-c).$$

It therefore suffices to bound $D_{n1}$ and $D_{n2}$ on $(-\infty, a_n]$. By Lemma 3.1, we may replace $C_n$ in the denominator by $C$. Hence, with large probability and up to a constant factor, (3.4) is bounded from above by

$$\int_0^{a_n} \frac{|C_n(y) - C(y)|}{C^2(y)} F_n^*(\mathrm{d}y) + \frac{1}{n} \int_0^{a_n} \frac{F_n^*(\mathrm{d}y)}{C^2(y)}.$$

Using (ii) of (A), it is easily seen that the expectation of the second term is bounded as $n \to \infty$. The expectation of the first term is less than or equal to

$$\int_0^{a_n} \frac{\mathbb{E}|C_{n-1}(y) - C(y)| + 1/n}{C^2(y)} F^*(\mathrm{d}y)$$

$$\leq \int_0^{a_n} \frac{(1/\sqrt{n-1})C^{1/2}(y) + 1/n}{C^2(y)} F^*(\mathrm{d}y) = O(1).$$

This proves (3.4). As for $D_{n1}$, we already mentioned that $D_{n1}$ converges to zero with probability one uniformly on each interval $[0, x_0]$ with $F(x_0) < 1$. Also, for each fixed $x$, we have

$$\mathbb{E}D_{n1}(x) = 0 \quad \text{and} \quad \text{Var}\, D_{n1}(x) \leq \frac{1}{n} \int_0^x \frac{F^*(\mathrm{d}y)}{C^2(y)}.$$

Moreover, by the construction of $a_n$ and the definitions of $F^*$ and $C$, we have

$$\text{Var}\, D_{n1}(a_n) = O(1).$$

We conclude that $D_{n1}(a_n) = O_{\mathbb{P}}(1)$. Applying standard tightness arguments for the (non-standardized) process $D_{n1}$ (see Billingsley (1968), page 128), we get that $D_{n1}$ is uniformly bounded on $[0, a_n]$. This, however, implies (3.3) and completes the proof of Lemma 3.2. □

Our next lemma implies that $S_{n1}$ is negligible.

**Lemma 3.3.** *Under the assumptions of Theorem 1.1, we have*

$$S_{n1} = o_{\mathbb{P}}(n^{-1/2}).$$

**Proof.** We first bound $B_n(x)$. Since, for $0 \leq x \leq \frac{1}{2}$, we have

$$-x - x^2 \leq -x - \frac{x^2}{2(1-x)} \leq \ln(1-x) \leq -x,$$

we obtain

$$-n^{-1} \int_0^{x-} \frac{F_n^*(\mathrm{d}y)}{C_n^2(y)} \leq -n \int_0^{x-} \frac{F_n^*(\mathrm{d}y)}{[nC_n(y) + 1]^2} \leq B_n(x) \leq 0. \tag{3.5}$$

Recall that $nC_n \geq 1$ on the support of $F_n^*$, so the above integrals are all well defined. We conclude from (3.5) that

$$|S_{n1}| \leq n^{-1} \int_0^\infty \frac{|\varphi(x)|}{C_n(x)} \mathrm{e}^{\Delta_n(x)} \int_0^{x-} \frac{F_n^*(\mathrm{d}y)}{C_n^2(y)} F_n^*(\mathrm{d}x).$$

Now, as in the proof of Lemma 3.2, for a given $\varepsilon > 0$, we may choose some small $c = c_\varepsilon$ and a sequence $a_n \to \infty$ such that $1 - F^*(a_n) = \frac{c}{n}$ and $\mathbb{P}(X_{n:n} \leq a_n) \geq 1 - \varepsilon$. Hence, on this event, $F_n^*$ has all of its mass on $[0, a_n]$ and integration with respect to $F_n^*$ may thus be restricted to $0 \leq x \leq a_n$. Furthermore, by Lemmas 3.1 and 3.2, the processes $C/C_n$ and $D_{n1} + D_{n2}$ remain stochastically bounded when restricted to the support of $F_n^*$ as $n \to \infty$. Together with $B_n(x) \leq 0$, it therefore suffices to show that

$$n^{-1/2} \int_0^{a_n} \frac{|\varphi(x)|\gamma(x)}{C(x)} \int_0^{x-} \frac{F_n^*(\mathrm{d}y)}{C^2(y)} F_n^*(\mathrm{d}x) = o_{\mathbb{P}}(1).$$

The expectation of the left-hand side is less than or equal to

$$n^{-1/2} \int_0^{a_n} \frac{|\varphi(x)|\gamma(x)}{C(x)} \int_0^{x-} \frac{F^*(\mathrm{d}y)}{C^2(y)} F^*(\mathrm{d}x).$$

For any fixed $x_0$ with $F(x_0) < 1$, the integral

$$\int_0^{x_0} \frac{|\varphi(x)|\gamma(x)}{C(x)} \int_0^{x-} \frac{F^*(\mathrm{d}y)}{C^2(y)} F^*(\mathrm{d}x)$$

is finite since $1 - F(y)$ is bounded away from zero there and, by assumption,

$$\int \frac{F^*(\mathrm{d}y)}{G^2(y)} = \alpha^{-1} \int \frac{\mathrm{d}F}{G} < \infty.$$

The same holds true for the integral

$$\int_{x_0}^{a_n} \frac{|\varphi(x)|\gamma(x)}{C(x)} \int_0^{x_0-} \frac{F^*(\mathrm{d}y)}{C^2(y)} F^*(\mathrm{d}x).$$

It remains to study

$$\int_{x_0}^{a_n} \frac{|\varphi(x)|\gamma(x)}{C(x)} \int_{x_0}^{x} \frac{F^*(\mathrm{d}y)}{C^2(y)} F^*(\mathrm{d}x).$$

This integral is bounded from above, however, by

$$G^{-1}(x_0) \int_{x_0}^{a_n} \frac{|\varphi(x)| F(\mathrm{d}x)}{1 - F(x)}.$$

Now apply Cauchy–Schwarz to get

$$\int_{x_0}^{a_n} \frac{|\varphi(x)| F(\mathrm{d}x)}{1 - F(x)} \leq \left[ \int_{x_0}^{\infty} \varphi^2(x) F(\mathrm{d}x) \right]^{1/2} \left[ \int_0^{a_n} \frac{F(\mathrm{d}x)}{[1 - F(x)]^2} \right]^{1/2}.$$

The second integral is less than or equal to $[1 - F(a_n)]^{-1}$. Since

$$\frac{c}{n} = 1 - F^*(a_n) = \alpha^{-1} \int_{a_n}^{\infty} G(z) F(\mathrm{d}z) \leq \alpha^{-1}(1 - F(a_n)),$$

the second square root is $O(n^{1/2})$. On the other hand, the first integral can be made arbitrarily small by choosing $x_0$ large enough. This concludes the proof of Lemma 3.3. □

Next, we study $S_{n2}$. For this, the following lemma will be crucial.

**Lemma 3.4.** *Under the assumptions of Theorem 1.1, we have*

$$I_n \equiv \int \frac{|\varphi|\gamma[C - C_n]^2}{C^2 C_n} \mathrm{d}F_n^* = o_{\mathbb{P}}(n^{-1/2}).$$

**Proof.** As in the proof of Lemma 3.3, we have $X_i \leq a_n$ for all $i = 1, \ldots, n$ with large probability. Similarly, consider a sequence $b_n > 0$ such that

$$F^*(b_n) = \frac{c}{n}, \quad c \text{ sufficiently small},$$

such that $X_i \geq b_n$ for $i = 1, 2, \ldots, n$ with large probability. In other words, up to an event of small probability, we may restrict integration in $I_n$ to the interval $[b_n, a_n]$. Also, in view of Lemma 3.1, the $C_n$ in the denominator may be replaced by $C$ (times a constant), with large probability. Hence, on a set with large probability, we have

$$I_n = n^{-1} \sum_{i=1}^{n} \frac{|\varphi(X_i)|\gamma(X_i)1_{\{b_n \leq X_i \leq a_n\}}}{C^2(X_i)C_n(X_i)}$$
$$\times \left[\frac{1}{n}\sum_{j\neq i}(1_{\{Y_j \leq X_i \leq X_j\}} - C(X_i)) + \frac{1}{n}(1 - C(X_i))\right]^2$$
$$\leq \frac{2}{n}\sum_{i=1}^{n}\frac{|\varphi(X_i)|\gamma(X_i)1_{\{\ldots\}}}{C^3(X_i)}\left\{\frac{1}{n}\sum_{j\neq i}(1_{\{\ldots\}} - C(X_i))\right\}^2$$
$$+ \frac{2}{n^3}\sum_{i=1}^{n}\frac{|\varphi(X_i)|\gamma(X_i)1_{\{\ldots\}}}{C^2(X_i)C_n(X_i)}. \tag{3.6}$$

The first sum has an expectation not exceeding

$$\frac{2}{n}\int_{b_n}^{a_n}\frac{|\varphi(x)|\gamma(x)}{C^2(x)}F^*(\mathrm{d}x).$$

Fix a small positive $x_1$ and, as in the previous proof, a large $x_0$. Assume, without loss of generality, that $b_n \leq x_1 \leq x_0 \leq a_n$. The integral

$$\int_{x_1}^{x_0}\frac{|\varphi|\gamma}{C^2}\mathrm{d}F^*$$

is finite. So, the middle part contributes an error $O(\frac{1}{n})$, which is smaller than desired. The upper part, $\int_{x_0}^{a_n}\ldots$, is dealt with as in the proof of Lemma 3.3, yielding a bound $o(n^{1/2})$ as $x_0$ gets large. The same holds true for the lower part. Finally, the second sum in (3.6) may be studied along the same lines as was the first, by starting with the inequality $nC_n(X_i) \geq 1$. The proof is thus complete. □

**Corollary 3.1.** *Under the assumptions of Theorem 1.1, we have*

$$S_{n2} = \int_0^\infty \frac{\varphi\gamma}{C}\mathrm{d}F_n^* + \int_0^\infty \frac{\varphi\gamma[C - C_n]}{C^2}\mathrm{d}F_n^* + o_\mathbb{P}(n^{-1/2}).$$

We shall come back to $S_{n2}$ later, but first proceed to $S_{n3}$. Recalling $D_{n1}$, we have

$$S_{n3} = -\int_0^\infty \frac{\varphi(x)\mathrm{e}^{\Delta_n(x)}}{C_n(x)}\int_0^{x-}\frac{F_n^*(\mathrm{d}y) - F^*(\mathrm{d}y)}{C(y)}F_n^*(\mathrm{d}x).$$

To get an expansion for $S_{n3}$, the next lemma will be crucial.

**Lemma 3.5.** *Under the assumptions of Theorem 1.1, we have*

$$II_n \equiv \int_0^\infty \frac{\varphi(x)\gamma(x)}{C_n(x)}[e^{\Delta_n(x)+\int_0^{x-}\frac{F^*(\mathrm{d}y)}{C(y)}} - 1]\int_0^{x-}\frac{F_n^*(\mathrm{d}y) - F^*(\mathrm{d}y)}{C(y)}F_n^*(\mathrm{d}x)$$
$$= o_\mathbb{P}(n^{-1/2}).$$

**Proof.** By Cauchy–Schwarz, on a set with large probability,

$$nII_n^2 \leq n^{1/2}\int_0^{a_n}\frac{|\varphi(x)|\gamma(x)}{C_n(x)}\left[\int_0^{x-}\frac{F_n^*(\mathrm{d}y) - F^*(\mathrm{d}y)}{C(y)}\right]^2 F_n^*(\mathrm{d}x) \quad (3.7)$$
$$\times n^{1/2}\int_0^{a_n}\frac{|\varphi(x)|\gamma(x)}{C_n(x)}[e^{\Delta_n(x)+\int_0^{x-}\frac{F^*(\mathrm{d}y)}{C(y)}} - 1]^2 F_n^*(\mathrm{d}y).$$

By Lemma 3.1, we may again replace $C_n$ by $C$. The expectation of the resulting first integral is then less than or equal to

$$n^{-1/2}\int_0^{a_n}\frac{|\varphi(x)|\gamma(x)}{C(x)}\int_0^{x-}\frac{F^*(\mathrm{d}y)}{C^2(y)}F^*(\mathrm{d}x),$$

which was already shown to be $o(1)$ in the proof of Lemma 3.3. It therefore remains to show that the second factor in (3.7) is also $o_\mathbb{P}(1)$. Putting

$$z_n(x) = \Delta_n(x) + \int_0^{x-}\frac{F^*(\mathrm{d}y)}{C(y)},$$

we may write

$$[e^{z_n(x)} - 1]^2 = z_n^2(x)e^{2\tilde{z}_n(x)},$$

where $\tilde{z}_n(x)$ is between zero and $z_n(x)$. Recalling that $\Delta_n(x)$ is between the first integral in $B_n(x)$ and $-\int_0^{x-}\frac{F^*(\mathrm{d}y)}{C(y)}$, and that $B_n(x) \leq 0$, we may infer that $\tilde{z}_n(x)$ is uniformly bounded from above in probability, on the support of $F_n^*$, as $n \to \infty$. Hence, it suffices to bound the term

$$n^{1/2}\int_0^{a_n}\frac{|\varphi(x)|\gamma(x)z_n^2(x)}{C(x)}F_n^*(\mathrm{d}x),$$

which, in turn, is less than or equal to

$$n^{1/2}\int_0^{a_n}\frac{|\varphi(x)|\gamma(x)}{C(x)}[B_n(x) + D_{n1}(x) + D_{n2}(x)]^2 F_n^*(\mathrm{d}x)$$
$$\leq n^{1/2}2\int_0^{a_n}\frac{|\varphi(x)|\gamma(x)}{C(x)}[B_n^2(x) + (D_{n1}(x) + D_{n2}(x))^2]F_n^*(\mathrm{d}x).$$

By (3.5),

$$\int_0^{a_n}\frac{|\varphi(x)|\gamma(x)B_n^2(x)}{C(x)}F_n^*(\mathrm{d}x) \leq n^{-2}\int_0^{a_n}\frac{|\varphi(x)|\gamma(x)}{C(x)}\left[\int_0^{x-}\frac{F_n^*(\mathrm{d}y)}{C_n^2(y)}\right]^2 F_n^*(\mathrm{d}x).$$

To bound the right-hand side, first replace $C_n$ with $C$. The squared term may then be viewed as a V-statistic, of which the leading term is a U-statistic. Its expectation is

$$\frac{n-1}{n}\left[\int_0^{x-}\frac{F^*(\mathrm{d}y)}{C^2(y)}\right]^2.$$

Thus, we have to show that

$$n^{-2}\int_0^{a_n}|\varphi(x)|\left[\int_0^{x-}\frac{F^*(\mathrm{d}y)}{C^2(y)}\right]^2 F(\mathrm{d}x) = o(n^{-1/2}).$$

This follows as in the proof of Lemma 3.3. Only the powers of $1 - F(x)$ are different. Finally, the error terms involving $D_{n1}$ and $D_{n2}$ may be dealt with similarly. The proof is thus complete. $\square$

**Lemma 3.6.** *Under the assumptions of Theorem 1.1, we have*

$$\int\frac{\varphi\gamma[C_n - C](x)}{CC_n(x)}\int_0^{x-}\frac{F_n^*(\mathrm{d}y) - F^*(\mathrm{d}y)}{C(y)}F_n^*(\mathrm{d}x) = o_{\mathbb{P}}(n^{-1/2}).$$

**Proof.** As above. $\square$

Lemmas 3.5 and 3.6 immediately imply the following corollary, which brings us to the desired representation of $S_{n3}$.

**Corollary 3.2.** *Under the assumptions of Theorem 1.1, we have the representation*

$$S_{n3} = -\int_0^\infty\frac{\varphi(x)\gamma(x)}{C(x)}\int_0^{x-}\frac{F_n^*(\mathrm{d}y) - F^*(\mathrm{d}y)}{C(y)}F_n^*(\mathrm{d}x) + o_{\mathbb{P}}(n^{-1/2}).$$

We now proceed to $S_{n4}$. As a first step to get the desired representation, we shall need the following lemma.

**Lemma 3.7.** *Under the assumptions of Theorem 1.1, $S_{n4}$ admits the expansion*

$$S_{n4} = \int_0^\infty\frac{\varphi(x)\mathrm{e}^{\Delta_n(x)}}{C_n(x)}\int_0^{x-}\frac{C_n(y) + 1/n - C(y)}{C^2(y)}F_n^*(\mathrm{d}y)F_n^*(\mathrm{d}x) \qquad (3.8)$$
$$+ o_{\mathbb{P}}(n^{-1/2}).$$

**Proof.** Recalling the definition of $D_{n2}$, the difference between $S_{n4}$ and the leading term in (3.8) becomes

$$-\int_0^\infty\frac{\varphi(x)\mathrm{e}^{\Delta_n(x)}}{C_n(x)}\int_0^{x-}\frac{[C_n(y) + 1/n - C(y)]^2}{C^2(y)[C_n(y) + 1/n]}F_n^*(\mathrm{d}y)F_n^*(\mathrm{d}x).$$

Taking Lemma 3.1 into account, the last expression is bounded in absolute values from above, with large probability, by

$$\int_0^\infty \frac{|\varphi(x)|e^{\Delta_n(x)}}{C(x)} \int_0^{x-} \frac{[C_n(y)+1/n-C(y)]^2}{C^3(y)} F_n^*(\mathrm{d}y) F_n^*(\mathrm{d}x).$$

In the proof of Lemma 3.5, we argued that

$$z_n(x) = \Delta_n(x) + \int_0^{x-} \frac{F^*(\mathrm{d}y)}{C(y)}$$

is uniformly bounded from above, on the support of $F_n^*$, as $n \to \infty$, with probability one. Consequently, we may replace $\exp(\Delta_n(x))$ by $\gamma(x)$. Also, we may wish to extend the integral only up to $a_n$. The expectation of the resulting term then does not exceed

$$\int_0^{a_n} |\varphi(x)| \int_0^{x-} \frac{(2/(n-1))C(y)+4/n^2}{C^3(y)} F^*(\mathrm{d}y) F(\mathrm{d}x).$$

These terms already appeared in the proof of Lemma 3.3 and were there shown to be $o(n^{-1/2})$. The proof is therefore complete. □

In the following, we shall omit the summand $\frac{1}{n}$ in (3.8), since its contribution is also negligible. The next lemma will enable us to replace $\exp(\Delta_n(x))$ by $\gamma(x)$.

**Lemma 3.8.** *Under the assumptions of Theorem 1.1, we have*

$$III_n \equiv \int_0^\infty \frac{\varphi(x)\gamma(x)}{C_n(x)} [e^{\Delta_n(x)+\int_0^{x-} \frac{F^*(\mathrm{d}y)}{C(y)}} - 1] \int_0^{x-} \frac{C_n - C}{C^2} \mathrm{d}F_n^* F_n^*(\mathrm{d}x)$$
$$= o_\mathbb{P}(n^{-1/2}).$$

**Proof.** Cauchy–Schwarz leads to a bound for $nIII_n^2$ similar to (3.7) for $nII_n^2$. The second factor is exactly the same as before. The first factor is dealt with in the following lemma. □

**Lemma 3.9.** *Under the assumptions of Theorem 1.1, we get*

$$\int_0^\infty \frac{|\varphi(x)|\gamma(x)}{C(x)} \left[\int_0^{x-} \frac{C_n(y)-C(y)}{C^2(y)} F_n^*(\mathrm{d}y)\right]^2 F_n^*(\mathrm{d}x) = o_\mathbb{P}(n^{-1/2}).$$

**Proof.** As in previous proofs, it is enough to extend the $x$-integral up to $a_n$. It is then easily checked that the expectation of the resulting term is less than or equal to (when $n > 3$)

$$\int_0^{a_n} |\varphi(x)| \int_0^{x-} \int_0^{x-} \frac{\frac{1}{(n-3)\alpha} G(y \wedge z)(1-F(y \vee z)) + 4/n^2}{C^2(y)C^2(z)} F^*(\mathrm{d}y) F^*(\mathrm{d}z) F(\mathrm{d}x),$$

where

$$y \wedge z = \min(y, z) \quad \text{and} \quad y \vee z = \max(y, z).$$

Neglecting $4n^{-2}$ for a moment, the rest equals

$$\frac{2}{(n-3)\alpha} \int_0^{a_n} |\varphi(x)| \int_0^{x-} \int_y^{x-} \frac{G(y)(1-F(z))}{C^2(y)C^2(z)} F^*(\mathrm{d}y) F^*(\mathrm{d}z) F(\mathrm{d}x)$$

$$= \frac{2}{(n-3)} \int_0^{a_n} |\varphi(x)| \int_0^{x-} \frac{1}{(1-F(y))^2} \int_y^{x-} \frac{1}{C(z)} F(\mathrm{d}z) F(\mathrm{d}y) F(\mathrm{d}x)$$

$$= \frac{2}{n-3} \int_0^{a_n} |\varphi(x)| \int_0^{x-} \frac{1}{C(z)} \int_0^z \frac{1}{(1-F(y))^2} F(\mathrm{d}y) F(\mathrm{d}z) F(\mathrm{d}x)$$

$$\leq \frac{2}{n-3} \int_0^{a_n} |\varphi(x)| \int_0^{x-} \frac{1}{C(z)(1-F(z))} F(\mathrm{d}z) F(\mathrm{d}x).$$

The last integral already appeared in the proof of Lemma 3.3 and was shown to be $o(n^{1/2})$. The error term $\frac{4}{n^2}$ yields similar bounds, so the proof is complete. □

In the final step, $C_n$ in the denominator of (3.8) may be replaced by $C$. The details are omitted. Hence, we arrive at the following representation of $S_{n4}$.

**Corollary 3.3.** *Under the assumptions of Theorem 1.1, we have*

$$S_{n4} = \int_0^\infty \frac{\varphi(x)\gamma(x)}{C(x)} \int_0^{x-} \frac{C_n(y) - C(y)}{C^2(y)} F_n^*(\mathrm{d}y) F_n^*(\mathrm{d}x)$$
$$+ o_\mathbb{P}(n^{-1/2}).$$

We are now ready to prove Theorem 1.1.

**Proof of Theorem 1.1.** Corollaries 3.1, 3.2 and 3.3 yield representations of the relevant terms $S_{n2}, S_{n3}$ and $S_{n4}$ as V-statistics. As is known (see Serfling (1980)), a V-statistic equals a U-statistic, up to an error $o_\mathbb{P}(n^{-1/2})$. If, in $S_{n2}, S_{n3}$ and $S_{n4}$, we replace $F_n^*$ by $F^*$, we come up with degenerate U-statistics which are all of the order $o_\mathbb{P}(n^{-1/2})$. Collecting the leading terms and applying Fubini's theorem then yields the proof of Theorem 1.1 with $\hat{F}_n$ instead of $F_n$.

As for $F_n$, apply the inequality

$$\left| \prod_{j=1}^{i-1} a_j - \prod_{j=1}^{i-1} b_j \right| \leq \sum_{j=1}^{i-1} |a_j - b_j|, \qquad |a_j|, |b_j| \leq 1,$$

to the weights of $F_n$ and $\hat{F}_n$. It follows that

$$\left| \int \varphi \, \mathrm{d}\hat{F}_n - \int \varphi \, \mathrm{d}F_n \right| \leq n^{-3} \sum_{i=1}^{n} \frac{|\varphi(X_{i:n})|}{C_n(X_{i:n})} \sum_{j=1}^{i-1} \frac{1}{C_n^2(X_{j:n})}$$

$$= n^{-1} \int \frac{|\varphi(x)|}{C_n(x)} \int_0^{x-} \frac{1}{C_n^2(y)} F_n^*(\mathrm{d}y) F_n^*(\mathrm{d}x).$$

Use Lemma 3.1 again to replace $C_n$ with $C$ and apply the SLLN (for U-statistics) to get that the last term is $O_{\mathbb{P}}(n^{-1})$. This completes the proof of Theorem 1.1. □

**Proof of Lemma 1.1.** The proof proceeds by induction on $i$. For $i=1$ and $j=1$, the left-hand side of (1.16) equals $1/nC_n^U(U_{11:n})$, as does the right-hand side. The assertion then follows from (1.14). For $i=1$ and $j=2$, the left-hand side of (1.16) equals, by (1.8) and (1.15),

$$\frac{(nC_n^U(U_{11:n}) - 1)/(nC_n^U(U_{11:n}))}{nC_n^U(U_{12:n})} = \frac{1}{nC_n^U(U_{11:n})} = \frac{1}{nC_n(X_{1:n})}.$$

For $i=1$ and $j>2$, the proof follows the same pattern, by repeated use of (1.15). Hence, the assertion of Lemma 1.1 holds true for $i=1$. Assuming that (1.16) holds true for all indices less than or equal to $i$, we obtain, by (1.7),

$$1 - F_n(X_{i:n}) = 1 - F_n^U(U_{id_i:n}).$$

This may be used as a starting point to show (1.16) for $i+1$. Again, make repeated use of (1.14) and (1.15). □

## Acknowledgement

The research of Jane-Ling Wang was supported by an NSF Grant NSF-04-06430.